\documentclass[reqno]{amsart}

\usepackage{amsfonts,amsmath,amssymb,amsrefs,graphicx,float}
\usepackage{enumerate}

\usepackage{tikz}
\usetikzlibrary{decorations.markings}
\usetikzlibrary{patterns}

\newcommand{\be}{\begin{equation*}}
\newcommand{\ee}{\end{equation*}}
\newcommand{\beq}{\begin{equation}}
\newcommand{\eeq}{\end{equation}}
\newcommand{\bs}{\begin{split}}
\newcommand{\esplit}{\end{split}}
\newcommand{\begincal}{\begin{eqnarray*}}
\newcommand{\fincal}{\end{eqnarray*}}

\newtheorem{thm}{Theorem}[section]
\newtheorem{lemma}{Lemma}[section]

\newtheorem{prop}{Proposition}[section]
\newtheorem{defi}{Definition}[section]

\newtheorem{rem}{Remark}[section]

\newcommand{\R}{{\mathbb R}}

\newcommand{\Sun}{{\mathbb{S}}^1}

\begin{document}

\title[An obstruction to the existence of prescribed curvature ]{An obstruction to the existence of immersed curves of prescribed curvature}
\author{Stephane Kirsch}
\address{Stephane Kirsch, UBC, PIMS , West Mall Annex, 1933 WEST MALL V6T 1Z2, Vancouver BC, CANADA}
\email{kirsch@math.ubc.ca}

\author{Paul Laurain }
\address{Paul Laurain, Ecole normale sup\'erieure de Lyon, D\'epartement de 
Math\'ematiques - UMPA, 46 all\'ee d'Italie, 69364 Lyon cedex 07,
France}
\email{Paul.Laurain@umpa.ens-lyon.fr}

\date{March 2009}

\begin{abstract}
We prove in this paper that there exists no immersed closed curve in $\R^2$ with curvature increasing in a given direction. This generalizes in the 1-dimensional setting a result of the first author for embedded hypersurfaces in $\R^{n+1}$.
\end{abstract}

\maketitle
\section{Introduction}
The problem of prescribing curvature of manifolds has been intensively studied in the last decades. In this paper we consider the problem of prescribing the mean curvature of a hypersurface in the Euclidean space $\R^n$. More precisely, considering a function $H:\R^n \to \R$, we can ask whether or not there exists a compact closed (embedded or immersed) hypersurface whose mean curvature at any point $x$ is given by $H(x)$. We are interested here in finding an obstruction to the existence of such a hypersurface, in the spirit of the Kazdan-Warner obstruction for scalar curvature on the 2-sphere. We look for conditions on $H$ ensuring that the problem of prescribed mean curvature has no solution. For example the first author proved in his thesis the following theorem:
\begin{thm}\label{teoemb}
Let $n \geq 2$ and $H: \R^n \to \R$ be a Lipschitz function such that there exists $\vec{e} \in \R^n$ verifying
$$
\langle \nabla H (x) , \vec{e} \rangle >0 \ \ \text{ for almost every } \ x \in \R^n,
$$
then there exists no $C^1$-closed hypersurface embedded in $\R^n$ of mean curvature $H$.
\end{thm}

In fact this theorem is a little bit more general since it still holds if we replace the mean curvature by any k-th mean curvature (the k-th symmetric function of the principal curvatures). Theorem \ref{teoemb} means that an embedded hypersurface can not have a mean curvature which is monotone in some direction. This obstruction is very similar to the obstruction found by Kazdan-Warner for the scalar curvature on the 2-sphere (see section 8 of \cite{KazdanWarner74} and also \cite {BourguignonEzin87}): there exists no metric on $S^2$ whose scalar curvature is monotone in a fixed direction of $\R^3$. An analogue of the Kazdan-Warner obstruction was recently proved by Ammann, Humbert and Ould Ahmedou in the mean curvature context (see \cite {AmmannHumbert08}): if one sees the mean curvature as a function on $S^n$ instead of the ambient space $\R^{n+1}$, there exists no conformal immersion of $S^n$ in $\R^{n+1}$ such that the mean curvature $H:S^n \to \R$ is monotone in one direction.\\

In this paper we investigate a generalization of theorem \ref{teoemb} in the immersed case. The immersed world is much richer and more complicated than the embedded world. For instance, Alexandrov proved (see \cite{Aleksandrov56}) that the only closed compact hypersurfaces of constant mean curvature embedded in $\R^n$ are the spheres. Though for many years it was conjectured that the embedding assumption could be removed, Wente constructed in the 80's (see \cite{Wente87}) immersed tori of constant mean curvature which are not embedded in any Euclidean space of dimension $n \geq 3$. We can then expect the extension of theorem \ref{teoemb} to be tricky (maybe even false) in the immersed case. Nevertheless the Alexandrov theorem holds for immersed spheres (see \cite{Hopf}) so we can expect some rigidity of spheres with respect to their mean curvature. In particular, the Alexandrov  theorem is true for immersed curved. This suggests to look first at an extension of theorem \ref{teoemb} in the 1-dimensional case. In this paper we prove the following:

\begin{thm} 
\label{t1} 
Let $H:\R^2 \to \R$ be a positive Lipschitz function such that there exists $\vec{e} \in \R^2$ verifying
$$
\left \langle\nabla H (x) , \vec{e} \right\rangle > 0 \ \ \text{ for almost every } x \in \R^2. 
$$
Then there exists no $C^1$-closed  curve immersed in $\R^2$ whose curvature at every point $x$ is $H(x)$.
\end{thm}
\begin{rem}
Without the positivity assumption on $H$ the result does not hold, whereas it still holds in the embedded case. Indeed, a simple counterexample is given by the lemniscate of Bernoulli.\\

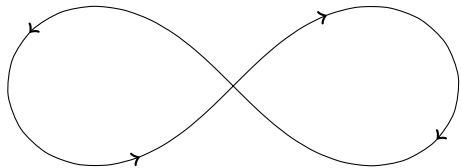
\begin{figure}[!h]
\centering
\begin{tikzpicture} [scale=1.5, arrow/.style={postaction={decorate}, decoration={markings,mark= at position 0.2 with {\arrow[line width=1pt]{>}}}, decoration={markings, mark= at position 0.4 with {\arrow[line width=1pt]{>}}}, decoration={markings, mark= at position 0.6 with {\arrow[line width=1pt]{>}}}, decoration={markings, mark= at position 0.8 with {\arrow[line width=1pt]{>}}}}]
\draw[domain=-3.141:3.141,samples=40,smooth,arrow]
plot[parametric] function{(2*(sin(t)))/(1+(cos(t))*(cos(t))),(sin(2*t))/(1+(cos(t))*(cos(t)))};
\end{tikzpicture}
\caption{Lemniscate of Bernouilli}
\label{Lemniscate}
\end{figure}
\vspace{0.15cm}
The cartesian equation of the lemniscate is 
$$
(x_1^2+x_2^2)^2 = a^2(x_1^2-x_2^2),
$$
where $a>0$ is a parameter. A straightforward computation gives the curvature $c(x_1)$ of the lemniscate:
$$
c(x_1) = sign(x_1) \frac{3}{a^2}\sqrt{\frac{a\sqrt{8x_1^2+a^2}-a^2}{2}},
$$
which is a $C^1$ increasing function.
The lemniscate of Bernoulli is thus a solution of the prescribed curvature problem where $H(x_1,x_2)=c(x_1)$, with
$$
\left\langle\nabla H(x) , e_1 \right\rangle = c'(x_1) > 0 \text{ for every } x \in \R^2.
$$
\end{rem}
\vspace{0.15cm}
Our proof relies on the fact that a self-intersection of a curve is very simple: it is a point. In higher dimension self-intersections of a hypersurface may be very complicated, so our proof seems hard to generalize, and in that case, the problem remains widely open.

The paper is organized as follows: in part 2, we give some definitions, notations, we state a preliminary lemma whose proof is given in the appendix, and we give a short sketch of
the proof of theorem \ref{t1}, and in part 3 we prove theorems \ref{teoemb} and \ref{t1}.

\section{Notations and preliminaries}
All curves are assumed to be at least in $C^2_p \cap C^0$ (i.e continuous, with a piecewise continuous first and second derivative). Moreover, thanks to the rescaling properties of our problem, we assume also that all curves have length $2\pi$ and are parametrized by arc length on $\Sun \simeq [0,2\pi[ $. We denote by $\Gamma$ a curve
parametrized by $\gamma: \Sun \rightarrow \R^2$. Let $\vec{t}(s)=\frac{d\gamma}{ds}(s)$ and $\vec{n}(s)$ be such that $(\gamma(s),\vec{t}(s),\vec{n}(s))$ is a positive orthonormal frame for every $s$ such that $\vec{t}(s)$ exists. Then for a $C^{2}_{p}$ curve, we can  write the equation of prescribed curvature on every $C^2$ piece of $\gamma$ as follows:
\beq
\label{presc}
\frac{d \vec{t}}{ds}(s)= H(\gamma(s)) \vec{n}(s).
\eeq
We define here the order of a self-intersection point:
\begin{defi}
Let $x$ be a point of $\Gamma$, we define the order of $x$, which we denote $ord(x)$, by
$$
ord(x) := \#\{ s \in \Sun \text{ such that } \gamma(s) = x \}.
$$
\end{defi}
It is obvious that $x$ is a self-intersection point if and only if $ord(x) \geq 2$.
\begin{rem}\label{intersections}
If $H:\R^2\rightarrow\R$ is a Lipschitz function, the Cauchy-Lipschitz theorem insures that every non-periodic $C^2$-solution of (\ref{presc}) has isolated self-intersection of finite order. Moreover, if two parts of the curve intersect tangentially then they must have opposite orientation.  
\end{rem}
Before going any further we need to define the index of a closed curve and some properties of it.
\begin{defi}
Let $\gamma \in C^0 ( \Sun, \R^2)$ and $a \in \R^2 \setminus \Gamma $, we define the index of $\gamma$ with respect to $a$ as the degree of the map from $S^1$ into $S^1$ defined by
\be
s \mapsto \frac{\gamma(s)-a}{\Vert \gamma(s)-a\Vert },
\ee
which we denote by $Ind_\gamma(a)$.
\end{defi}
\begin{rem}
For a given curve $\Gamma$ the index is defined up to a choice of orientation.
\end{rem}
It is well known that the map $x\mapsto Ind_\gamma(x)$ is constant on each connected component of $\R^2 \setminus \Gamma$. Moreover, thanks to the following proposition, we can quantify the variation of index passing from a component to another.\\
\begin{prop}
\label{saut}
Let $\gamma \in C^{1}_p ( \Sun, \R^2)\cap C^0( \Sun, \R^2)$ and $a \in \R^2 \setminus \Gamma$. Let $p\in \Gamma$ be a regular point; by regular point we mean that $p$ is not a self-intersection
point of $\Gamma$ and $\gamma$ is $C^1$ in a neighborhood of $\gamma^{-1}(p)$; then for $r$ sufficiently small, $B(p,r)\setminus \Gamma = U \cup V$, where $U$ and $V$ are two connected open subsets of $\R^2$ such that $U \cap V = \emptyset$ and $\vec{n}$ points into $U$. Moreover $\forall (u,v) \in U\times V$ we have
\be
Ind_{\gamma}(u)-Ind_{\gamma}(v) =1.
\ee   
\end{prop}
A proof  can be  found  in Lemma 9.2.5 of \cite{BergerGostiaux88}; moreover it is obviously false if $p$ is a point of self-intersection as shown by a circle parametrized twice.\\
We are now able to give the following lemma, which is a generalization of the Stokes theorem:
\begin{lemma}
\label{ldiv}
Let $\Gamma \subset \R^2$ be a continuous  $C^1_p$-closed curve with isolated self-intersections of finite order, and $\vec{V}$ be a Lipschitz vector field on $\R^2$, then we have the following formula:
\beq \label{eq1}
\int_{\R^2} div(\vec{V}) Ind_\gamma(x) dx = -\oint_\Gamma \left\langle\vec{V},\vec{n}\right\rangle d\sigma,
\eeq
where $d\sigma$ is the restriction of the Lebesgue measure to $\Gamma$ and $\vec{n}$ the unit normal vector field on $\Gamma$.
\end{lemma}
The proof of this lemma is postponed to the appendix.
\begin{rem}
There are two possible choices for the orientation of $\Gamma$, but (\ref{eq1}) obviously holds in both cases, because  turning $\vec{n}$ into $-\vec{n}$ turns $Ind_\gamma(x)$ into $-Ind_\gamma(x)$.  
\end{rem}
\begin{rem}
If $\Gamma$ has no self-intersection, we have the classical Stokes formula.
Indeed, let $\vec{n}$ be the outward unit normal vector, then
$Ind_\gamma(x)$ is $-1$ if $x$ is in the area enclosed by the curve, which we denote $\Omega$, and $0$
otherwise, which gives
$$
-\int_\Omega div(\vec{V})=-\int_{\partial \Omega}\left\langle\vec{V} , \vec{n}\right\rangle.
$$
\end{rem}

Formula (\ref{eq1}) enables us to prove theorem \ref{t1} in the case where the index is everywhere non negative (resp. everywhere non
positive). 

Let $\Gamma$ be a $C^1$-closed  curve of curvature $H$, and let $\vec{V}= H \vec{e}$ be such that $\left\langle\nabla H , \vec{e}\right\rangle>0$, then we get
\beq
\label{div}
\begin{split}
\int_{\R^2} \left\langle\nabla H , \vec{e}\right\rangle Ind_\gamma(x) dx &= -\oint_\Gamma \left\langle\vec{e},H\vec{n}\right\rangle d\sigma\\
&= -\left\langle\vec{e},\oint_\Gamma \frac{d \vec{t}}{ds} d\sigma \right\rangle\\ 
&=0,
\end{split}
\eeq
leading to $Ind_\gamma=0$ almost everywhere in $\R^2$, which is impossible.\\

\begin{figure}[!h]
\centering 
\begin{tikzpicture} [scale=2, arrow/.style={postaction={decorate}, decoration={markings,mark= at position 0.15 with {\arrow[line width=1pt]{<}}}, decoration={markings, mark= at position 0.35 with {\arrow[line width=1pt]{<}}}, decoration={markings, mark= at position 0.55 with {\arrow[line width=1pt]{<}}}, decoration={markings, mark= at position 0.75 with {\arrow[line width=1pt]{<}}}, decoration={markings, mark= at position 0.95 with {\arrow[line width=1pt]{<}}}}]
\draw[arrow] (1,0)
\foreach \x in {0,1,2,3,4} 
  {..controls (\x*72+10 :1.25) and (\x*72+7 :1.75) .. (\x*72 :1.75)..controls (\x*72-7 :1.75) and (\x*72-10 :1.25)  .. (\x*72 :1)..controls (\x*72+24 :0.7) and (\x*72+48 :0.7) .. (\x*72+72 :1)};
\foreach \i in {0,1,2,3,4}
\draw (\i*72 : 1.5) node{1};
\draw (0,0) node{-1};
\draw[line width=1.4pt] (216 :1)..controls (226 :1.25) and (223 :1.75) .. (216 :1.75)..controls (209 :1.75) and (206 :1.25)  .. (216 :1);
\end{tikzpicture}
\caption{Flower, a curve whose index has different signs}
\label{fleur}
\end{figure}
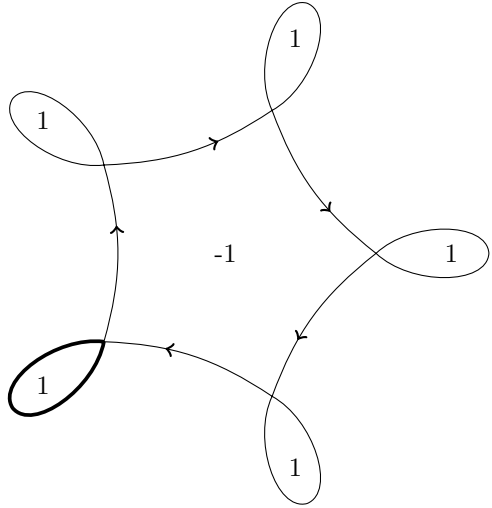

In order to prove the general case, we consider some closed subpart of the curve whose own index is non negative almost everywhere in $\R^2$. The idea is to cut off a loop of  the curve, as the bold one in figure \ref{fleur}. Of course this subcurve may not be regular and then 
$$ \oint_\Gamma \frac{d \vec{t}}{ds} d\sigma\not=0,$$ 
but to get a contradiction we just need to have the right-hand side of (\ref{div}) negative. For that purpose the singularities of our loop have to be ``well oriented'' (left side on figure \ref{fleur}).
The following section is devoted to the proof of theorem \ref{teoemb} and theorem \ref{t1} that we have just sketched.
 
\newpage

\section{Proof of the theorems}
\subsection{Proof of theorem \ref{teoemb}} 
Suppose by contradiction that there exists a closed compact $C^1$-hypersurface $S$ embedded in $\R^n$ whose mean curvature is $H$. Then using the first variation formula with a constant vector field we easily find the following identity:
\begin{equation}\label{mean}
\int_S H(x)\vec{n}(x)dS = 0,
\end{equation}
where $\vec{n}$ is the unit normal vector field on $S$, given by $H$, and $dS = dH^{n-1}|_S$. In the case of the k-th mean curvature the identity (\ref{mean}) still holds (see \cite{Chen71}). Since $S$ is embedded we can use the generalized Jordan-curve theorem (see \cite{Lima}) to find $\Omega \subset \R^n$ a bounded open set such that $S = \partial\Omega$. We have then, using the Stokes formula,
$$
\int_S H\vec{n}dS = \int_\Omega \nabla H(x)dx = 0,
$$
which contradicts the assumption $\langle \nabla H(x), \vec{e}\rangle >0$, and thus proves theorem \ref{teoemb}.\hfill$\square$
\subsection{Proof of theorem \ref{t1}}
We prove the theorem by contradiction. We assume that there exists a closed $C^2$-curve $\Gamma$ of curvature $H$. Without loss of generality, we may assume that $\vec{e}=(1,0)$. \\

We first give some definitions when the curve $\Gamma$ we consider is only in $C^0 \cap C_p^2$.\\

We say that $x$ is a {\bf singularity} with respect to  $s$ if $x=\gamma(s)$ and
$$
\Delta \vec{t} (s) :=\lim_{u\rightarrow 0^+} \vec{t}(s-u) - \vec{t}(s+u) \not= \vec{0}.
$$
Moreover we say that the singularity is {\bf well-oriented} if 
$$
\left\{
\begin{array}{ccc}
\displaystyle\left\langle\lim_{u\rightarrow 0^+} \vec{t}(s-u),\vec{e}\right\rangle & > & 0, \\
\displaystyle\left\langle\lim_{u\rightarrow 0^+} \vec{t}(s+u),\vec{e}\right\rangle & < & 0.
\end{array}
\right.
$$
Thus, for a well-oriented singularity, we have $\langle \Delta \vec{t} (s), \vec{e}\rangle > 0$.\\ 

As announced at the end of the previous section, in order to get a contradiction, we need to find $\Gamma' \subset \Gamma$ such that
$$ Ind_{\gamma'} \geq 0\hbox{ and } \left\langle\vec{e},\oint_{\Gamma'} \frac{d \vec{t}}{ds} d\sigma \right\rangle >0.$$
For that purpose, we construct by induction a decreasing sequence of closed subcurves $\Gamma_n$ of $\Gamma$, such that all the singularities of $\Gamma_n$ are well-oriented.\\

Starting  from $\Gamma_0=\Gamma$, we obtain $\Gamma_{n+1}$ from $\Gamma_{n}$ as follows : Consider$$\Omega_n :=\{x \in \R^2 \backslash \Gamma_n \hbox{ s.t. } Ind_{\gamma_n}(x)<0 \}.$$ If $\Omega_n=\emptyset$ then $\Gamma_{n+1}=\Gamma_n$, otherwise $\Omega_n$ is bounded and we cut $\Gamma_n$ at the ``left extremity" of $\Omega_n$. In order to do so, we proceed in several steps:\\

{\bf Step $1$: isolate the ``left extremity" of $\Omega_n$}\\

Let $x=(x_1,x_2)\in\overline{\Omega_n}$ be such that $x_1$ is minimal. We claim that $x$ is a point of self-intersection of $\Gamma_n$.
Otherwise it is either a regular point, or a singularity created at a previous step of the algorithm which is not a self-intersection point.\\
\begin{tabular}{ll}
\begin{minipage}{10cm}
Suppose it is a regular point, then thanks to remark \ref{intersections}, there is no self-intersection of $\Gamma$, and consequently of $\Gamma_n$, in a small ball centered at $x$. Moreover, the tangent line to $\Gamma_n$ at $x$ has to be vertical and  the curve lies locally on the right-hand side of the tangent line, because of the minimality of the abscissa of $x$. Now, because of the orientation imposed by $H>0$ and proposition \ref{saut}, the index has to decrease when we cross $\Gamma_n$ at the point $x$ from the right to the left, which contradicts again the minimality of the abscissa of $x$. Hence $x$ can not be a regular point. \\
\end{minipage}
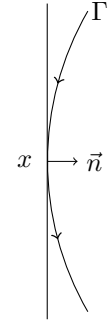
\begin{minipage}{2cm}
\begin{center}
\begin{tikzpicture}[arrow/.style={postaction={decorate}, decoration={markings,mark= at position 0.25 with {\arrow[line width=0.6pt]{>}}}, decoration={markings, mark= at position 0.75 with {\arrow[line width=0.6pt]{>}}}}]
\draw (0,-0.1) -- (0,4.1);
\draw[arrow] (0.54,4) arc (150:210:40mm);
\draw (0.7,4) node {$\Gamma$};
\draw [->,very thin] (0,2) -- (0.4,2) node [anchor = west] {$\vec{n}$};
\draw (-0.3,2) node {$x$};
\end{tikzpicture}

\end{center}
\end{minipage}
\end{tabular}\\
\begin{tabular}{ll}
\begin{minipage}{10cm}
Now, suppose $x$ is a singularity created at a previous step of the algorithm and not a self-intersection point. Since, as we shall see, the algorithm only creates well-oriented singularities---which means that the curve arrives transversally from the left (i.e $\{(x,y) \hbox{ s.t. } x\leq x_1\}$) on the vertical line passing through $x$ and goes back transversally to the left---it is easy to find a point $\hat{x}$ whose abscissa is smaller than the abscissa of $x$ and verifiying $Ind_{\gamma_n}(\hat{x})<0$, which is a contradiction and proves the claim.\\
\end{minipage}
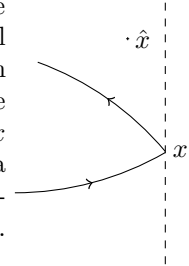
\begin{minipage}{2cm}
\begin{center}
\begin{tikzpicture}[arrow/.style={postaction={decorate}, decoration={markings,mark= at position 0.5 with {\arrow[line width=0.6pt]{>}}}}]
\draw[arrow] (0,0) arc (40:70:40mm);
\draw[arrow] (-2,-0.55) arc (270:300:40mm);
\draw[dashed] (0,-1.5) -- (0,2);
\draw (0.2,0) node {$x$};
\draw (-0.52,1.5) -- (-0.48,1.5);
\draw (-0.5,1.48) -- (-0.5,1.52);
\draw (-0.3,1.5) node {$\hat{x}$};
\end{tikzpicture}

\end{center}
\end{minipage}
\end{tabular}
\newpage
{\bf Step $2$: picture of the curve around $x$} \\

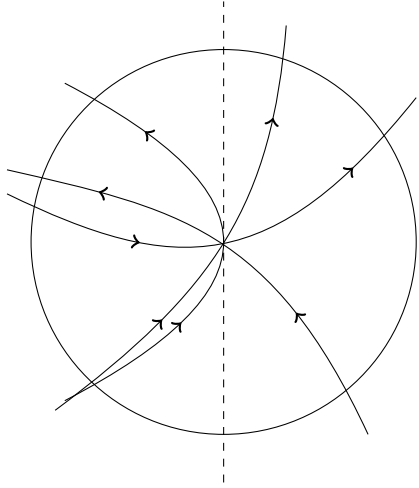
\begin{figure}[!h]
\centering
\begin{tikzpicture}[scale=0.64,arrow/.style={postaction={decorate}, decoration={markings, mark= at position .3  with {\arrow[line width=1pt]{>}}}, decoration={markings, mark= at position .8  with {\arrow[line width=1pt]{>}}}}]

\draw[dashed] (0,-5) to (0,5);
\draw (0,0) circle (4cm);
\draw[arrow] (-3.3,-3.3) .. controls (1.1,-1) and (1.1,1) .. (-3.3,3.3);
\draw[arrow] (3,-4).. controls (1,0.3)  and (-1,0.7)..(-4.5,1.5);
\draw[arrow]  (-4.5,1) .. controls (-1,-0.7) and (1,-0.75).. (4,3);
\draw[arrow] (-3.5,-3.5) .. controls (-0.15,-1) and (1,0.8).. (1.3,4.5);

\end{tikzpicture}  
\caption{Picture of the curve around $x$}
\label{sector}
\end{figure}

By remark \ref{intersections}, there exists a neighborhood of $x$ such that $\Gamma_n\setminus \{ x\}$ has no self-intersection. Then for $r>0$ small enough, the different parts of the curve divide the disk $D(x,r)$ into several sectors (see figure \ref{sector}).

\begin{figure}[!h]
\centering
\begin{tikzpicture}[scale=0.64, arrow/.style={postaction={decorate}, decoration={markings, mark= at position .3  with {\arrow[line width=1pt]{>}}}, decoration={markings, mark= at position .8  with {\arrow[line width=1pt]{>}}}}]

\draw[dashed] (0,-5) to (0,5);
\draw (0,0) circle (4cm);
\draw[arrow] (-3.3,-3.3) .. controls (1.1,-1) and (1.1,1) .. (-3.3,3.3);
\draw (3,-4).. controls (1,0.3)  and (-1,0.7)..(-4.5,1.5);
\draw  (-4.5,1) .. controls (-1,-0.7) and (1,-0.75).. (4,3);
\draw[arrow] (-3.5,-3.5) .. controls (-0.15,-1) and (1,0.8).. (1.3,4.5);
\draw(0,0) circle (4cm);
\draw[->,line width=1pt,xshift=2.4cm,yshift=1.3cm, rotate=44] (0,0) --  node[left, scale=1.2] {$\vec{n}$} (0,1);
\draw[->,line width=1pt,xshift=2.4cm,yshift=1.3cm, rotate=44] (0,0) -- node[above, scale=1.2] {$\vec{t}$}(1,0);
\draw[->,line width=1pt,xshift=2.1cm,yshift=-2.3cm, rotate=-55] (0,0) --  node[right, scale=1.2] {$\vec{n}$} (0,-1);
\draw[->,line width=1pt,xshift=2.1cm,yshift=-2.3cm, rotate=-55] (0,0) -- node[left, scale=1.2] {$\vec{t}$}(-1,0);
\filldraw [black, pattern=north east lines] (0,-0.04).. controls (1,0.16) and (2,0.8) .. (3.325,2.225).. controls (4,1.5)and (4.7,-1.25) ..(2.55,-3.08).. controls  (1.995,-2.1) and (1.095,-0.8) .. (0,-0.04) ;
\end{tikzpicture}
\caption{Sector of least index}
\label{cam2}
\end{figure}
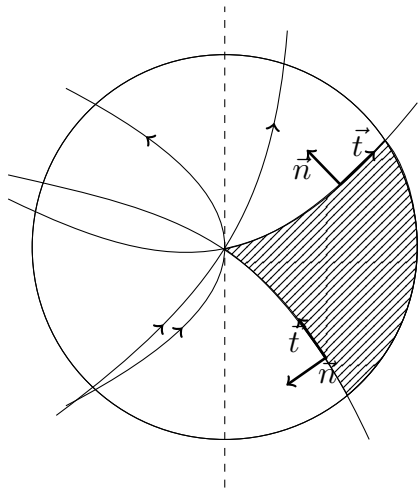

Let $(S_i)$ be an indexation of this partition of $D(x,r)\backslash\Gamma_n$. Let $S_0$ be the sector with the smallest index. First, we remark that this sector lies on the right-hand side of $x$, because its index is negative. This ensures that it cannot be edged by a singularity, because this one would not be well-oriented. Since the index has to increase crossing out of $S_0$, the normal vector, imposed by the orientation and $H>0$, points outward of the sector, as shown by figure \ref{cam2}. 

{\bf Step $3$:  cut the curve at $x$}\\

We choose an arc length parametrization $\gamma_n : [0,2\pi[ \rightarrow \R^2$ of $\Gamma_n$, such that $\gamma_n(0)\not=x$. Then there exist $\epsilon_1 >0$, $\epsilon_2>0$ and $s_{1}<s_{2}\in[0,2\pi[ $ such that the boundary of the sector of least index (where we do not consider the boundary of $D(x,r)$) is parametrized by two curves: $\gamma_n (]s_{1}-\epsilon_1, s_{1}])$ and $\gamma_n ([s_{2}, s_{2}+ \epsilon_2])$. Thanks to the fact that $x=\gamma_{n}(s_{1})=\gamma_{n}(s_{2})$ is not a singularity with respect to $s_1$ or $s_2$, we have $s_{1}\not= s_{2}$. We set  $\Gamma_{n+1} =\gamma_{n}(]s_1,s_2])$,  and we easily check that the set of singularities of $\Gamma_{n+1}$ is the singularity we created at $x$, which is well-oriented by construction, plus a subset of the singularities of $\Gamma_n$, which are all well-oriented by induction hypothesis. To end the algorithm, we rescale in order to have $\Gamma_{n+1}$ of length $2\pi$.\\

{\bf End of the proof:}\\

We have a decreasing sequence of curves $\Gamma_n$. Let $C_{n}$ be the number of connected components of $\R^{2} \setminus \Gamma_n$. $C_0$ is finite because $\Gamma$ has a finite number of self-intersection points  which are all of finite order. Moreover $C_n$ is decreasing, since at every step we cut off a part of the curve which encloses at least one component. But $\Gamma_{n}$ is never empty, because at each step we keep a part of the curve,  so the cutting process has to stop. This means that for $n$ large enough $Ind_{\gamma_{n}}$ is nonnegative on $\R^2\setminus\Gamma_{n}$, with $\Gamma_n \neq \emptyset$. For such a curve we get
\be
\left\langle\vec{e},\oint_{\Gamma_{n}} \frac{d \vec{t}}{ds} d\sigma \right\rangle = \sum_{i=1}^{k} \langle \Delta \vec{t}(s_i),\vec{e}\rangle >0,
\ee
where $\gamma_{n}(s_{i})$ are the singularity of $\Gamma_{n}$. As already remarked at the beginning of this section, this ends the proof.
\hfill$\square$
\appendix
\section{Proof of Lemma \ref{ldiv}}
Here we give a proof of  lemma \ref{ldiv} using elementary tools, though this can be proved shortly using degree theory, see for example \cite{fonseca}.\\
 
We denote $\Omega := \mathbb{R}^2 \backslash \Gamma$, which is an open subset of $\mathbb{R}^2$, and $\left(\Omega_i\right)_{i=0..n}$ its connected components, where
$\Omega_0$ is the unbounded one. We remark that on every $\Omega_i$, the index function is constant, because it changes only when we cross $\Gamma$, so we can define
$Ind_\gamma(\Omega_i)$ to be $Ind_\gamma(x)$ for every $x \in \Omega_i$. Hence 
\begin{eqnarray*}
\int_{\mathbb{R}^2}div(\vec{V})Ind_\gamma(x)dx & = & \sum_{i=0}^{n}\int_{\Omega_i}div(\vec{V}) Ind_\gamma(x)dx, \\
& = & \sum_{i=0}^{n} Ind_\gamma(\Omega_i)\int_{\Omega_i} div(\vec{V})dx, \\
& = & \sum_{i=0}^{n} Ind_\gamma(\Omega_i)\int_{\partial \Omega_i} \left\langle\vec{V} , \vec{n}_i\right\rangle dx, 
\end{eqnarray*}
where $\vec{n}_i$ is the outward normal vector to $\Omega_i$. The last equality is a consequence of the standard Stokes's theorem.

We can remark that $\Omega_0$ is not bounded and has infinite Lebesgue measure, but $Ind_\gamma(\Omega_0)=0$, so the term for $i=0$ in the right-hand side of the equation is $0$.

The boundary of each $\Omega_i$ is a $C^1_p$ curve that we can decompose as follows
$$
\partial \Omega_i = \bigcup_{j=0}^{m_i} \Gamma_{i,j},
$$
where each $\Gamma_{i,j}$ is a $C^1_p$ curve, and the self-intersections of $\Gamma$ belonging to $\partial \Omega_i$ are the points where we switch from a
$\Gamma_{i,j}$ to another.\\

Coming back to our computations we get
\begin{equation}\label{eq2}
\int_{\mathbb{R}^2}div(\vec{V})Ind_\gamma(x)dx = \sum_{i=0}^{n}\sum_{j=0}^{m_i} Ind_\gamma(\Omega_i)\int_{\Gamma_{i,j}} \left\langle\vec{V} , \vec{n}_{i,j} \right\rangle dx,
\end{equation}
where $\vec{n}_{i,j}$ is the outward normal vector to $\Omega_i$ on $\Gamma_{i,j}$.\\

We remark that in (\ref{eq2}), each $\Gamma_{i,j}$ appears exactly twice because it is on the boundary of exactly two $\Omega_i$. Thus for each couple $(i,j)$ there exists a
unique couple $(i',j')$, with $i\neq i'$ such that $\Gamma_{i,j}=\Gamma_{i',j'}$.\\

\begin{figure}[!h]
\centering
\begin{tikzpicture}
\draw (0mm,0mm) arc (-30:30:40mm);
\draw[->,very thin] (0.54,2) -- (0.85,2) node [anchor = west]{$\vec{n}_{i,j}=-\vec{n}_{i',j'}$};
\draw (2,3) node {$\Omega'_i$};
\draw (-1,2) node {$\Omega_i$};
\draw (1.1,0) node {$\Gamma_{i,j}=\Gamma_{i',j'}$}; 
\end{tikzpicture}
\caption{$\Gamma$ between two connected components}
\label{annexe} 
\end{figure}
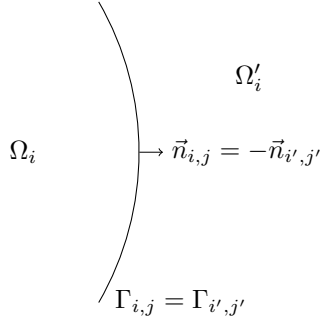

We can now change the indices of the sum and take a sum over $k$ in order to count each $\Gamma_{i,j}$ only once.
$$
\int_{\mathbb{R}^2}div(\vec{V})Ind_\gamma(x)dx = \sum_k \int_{\Gamma_k}
Ind_\gamma(\Omega_i) \left\langle\vec{V} , \vec{n}_{i,j} \right\rangle +
Ind_\gamma(\Omega_{i'})\left\langle\vec{V} , \vec{n}_{i',j'} \right\rangle .
$$
Moreover, we have in the previous formula that $\vec{n}_{i,j} =
-\vec{n}_{i',j'}$ because the outward normal for the open on one side of
$\Gamma_k$ is the inward normal for the open on the other side. We get then
$$
\int_{\mathbb{R}^2}div(\vec{V})Ind_\gamma(x)dx = \sum_k \int_{\Gamma_k}
\left(Ind_\gamma(\Omega_i)-Ind_\gamma(\Omega_{i'})\right) \left\langle\vec{V} , \vec{n}_{i,j}\right\rangle.
$$
Finally, by definition of $\vec{n}$, we have that $\vec{n}$ points in the
direction of increasing index. Moreover, since we assume only pointwise
self-intersections for $\gamma$, we have
$Ind_\gamma(\Omega_i)-Ind_\gamma(\Omega_{i'})=\pm 1$. More precisely we have
$$
\left\{
\begin{array}{ccccccc}
Ind_\gamma(\Omega_i)-Ind_\gamma(\Omega_{i'}) & = & -1 & \text{ if } & \vec{n}_{i,j} & = & \vec{n},\\
Ind_\gamma(\Omega_i)-Ind_\gamma(\Omega_{i'}) & = & 1 & \text{ if } & \vec{n}_{i,j} & = &-\vec{n}.
\end{array}
\right.
$$
Consequently we get
$$
(Ind_\gamma(\Omega_i)-Ind_\gamma(\Omega_{i'})) \vec{n}_{i,j}=-\vec{n}.
$$
Since $(\Gamma_k)_k$ is a partition of $\Gamma$, we obtain
$$
\int_{\mathbb{R}^2}div(\vec{V})Ind_\gamma(x)dx = -\sum_k \int_{\Gamma_k}
\left\langle\vec{V} , \vec{n} \right\rangle = -\oint_\Gamma \left\langle\vec{V} , \vec{n} \right\rangle.
$$
This ends the proof.$\hfill \Box$\\

\medskip {\bf Acknowledgements} : This work started while the second author was staying at PIMS, Vancouver. He wishes to thank all the members of PIMS for their great hospitality. The two authors would like to thank Olivier Druet for interesting and stimulating discussions.

\bibliographystyle{plain}
\bibliography{KirschLaurain09}

\end{document}